\documentclass[12pt]{article}
\usepackage{amsmath}
\usepackage{amssymb}
\usepackage{graphicx}
\usepackage[cp1251]{inputenc}
\usepackage[russian]{babel}
\newcommand{\il}[2]{\int\limits_{#1}^{#2}}

\newcommand{\ph}{\phantom{a}}
\newcommand{\phh}{\phantom{aaa}}

\newcommand{\sist}[2]{\left\{
\begin{array}{l}
{#1}\\
\ph\\
{#2}
\end{array}
\right.}

\begin{document}

\vskip 20pt

MSC 34C10

\vskip 20pt

\centerline{\bf Global solvability criterion for}
 \centerline{\bf matrix Riccati equations}

\vskip 20 pt

\centerline{\bf G. A. Grigorian}
\centerline{\it Institute  of Mathematics NAS of Armenia}
\centerline{\it E -mail: mathphys2@instmath.sci.am}
\vskip 20 pt

\noindent
Abstract. A new approach is used to obtain a global solvability criterion for matrix Riccati equations. It is shown that the obtained result is an extension of a result derived from a comparison theorem for matrix Riccati equations. Two corollaries   were drawn  from the obtained result as well.

\vskip 20 pt

Key words: matrix Riccati equations,  nonnegative (positive) definiteness, Hermitian matrices, The Liouville's formula.

\vskip 20 pt

{\bf 1. Introduction}.
Let $P(t), \ph Q(t), \ph R(t)$ and $S(t)$ be complex valued locally integrable   matrix  functions of dimension $n \times n$ on $[t_0,+\infty)$ and let $P(t)$ be Hermitian,  i.e.  $P(t) = P^*(t), \ph t\ge t_0$, where $*$ denotes the conjugation sign. Consider the  matrix Riccati equation
$$
Y' + Y P(t) Y + Q(t) y + Y R(t) - S(t) = 0, \phh t \ge t_0. \eqno (1.1)
$$
 By a solution  of this equation on $[t_1,t_2) \subset [t_0,+\infty)$ we mean an absolutely continuous matrix functions $Y = Y(t)$ of dimension $n\times n$ on $[t_1,t_2)$ satisfying (1.1) almost everywhere on $[t_1,t_2)$. It follows from the general theory of normal systems of differential equations that for every $t_1 \ge t_0)$ and for every initial value $Y_0$ there exists $t_2 > t_1$ such that Eq. (1.1) has a solution $Y(t)$ on $[t_1,t_2)$ with $Y(t_0) = Y_0$. Main interest from the point of view of qualitative theory of differential equations represents the case when $t_2 = +\infty$. In this case one says  that Eq. (1.1) is global solvable or has a global solution. An important problem is to find explicit conditions on the coefficients of Eq. (1.1), providing existence of global solutions of Eq. (1.1). This problem has been studied by many authors, and many works are devoted to it (see [2], and cited works therein).
In the scalar (real, complex and quaternionic) case the global solvability problem  for Eq. (1.1) was studies in the papers [6--8. 10] and the obtained results have been found some applications in [6--10].
In the general case the global solvability problem for Eq. (1.1) appear, in particular, in applied mathematics, notably in variational theory, optimal control and filtering, dynamic programming and differential games, invariant embedding  and scattering processes.
Knobloch and Pohl [13] have obtained  a global existence criterion for Eq. (1.1) based on a special maximum principle for second-order  linear differential equations. In [3] the authors use a Lyapunov  type function for obtaining a global existence result for Eq. (1.1). Reid [14] used an element wise comparison method in order to obtain a global existence theorem for nonsymmetric matrix Riccati  differential equations. This method has been used, modified and extended in the works [1, 5, 11, 12]. Among these results of this direction notice the following  (we represent this result here in a form convenient for comparison with the obtained result of this paper).

 For any Hermitian matrices $H_1$ and $H_1$ (of the same dimension) we denote by $H_1 \ge H_2 \ph (>0)$ the nonnegative (positive) definiteness of $H_1 - H_2$, by the symbol $H_1 \le 0 \ph (<0)$ we denote the nonnegative (positive) definiteness of $-H_1$.

{\bf Theorem 1.1. ([2, Theorem 3.5])}. {\it If $P(t) \ge 0, \ph S(t) \ge 0, \ph R(t) = Q^*(t), \phh t \ge t_0$, then the unique solution $Y(t)$ of Eq. (1.1) with piecewise continuous and locally integrable  coefficients exists for $t \ge t_0$ with
$$
0\le Y(t) \le \widetilde{Y}(t), \phh t \ge t_0.
$$
where $\widetilde{Y}(t)$ is the solution of the linear equation
$$
Y' + A^*(t) Y + Y A(t) - S(t) = 0
$$
with $\widetilde{Y}(t_0) = Y_0 \ge 0$.}

In this paper we extend  (generalize) this result up to the estimate $Y(t) \le \widetilde{Y}(t), \ph t \ge t_0.$

\vskip 5pt

{\bf 2. Auxiliary propositions}. The next four lemmas play crucial role in the proof of the main result. Let $M_k \equiv (m_{ij}^l)_{i,j=1}^n, \ph l=1,2$ be complex-valued matrices.

{\bf Lemma 2.1.} {\it The equality
$$
tr(M_1 M_1) = tr (M_2 M_1)
$$
is valid.}

Proof.  We have $tr (M_1 M_2) = \sum\limits_{j=1}^n(\sum\limits_{k=1}^n m_{jk}^1 m_{kj}^2) = \sum\limits_{k=1}^n(\sum\limits_{j=1}^n m_{jk}^1 m_{kj}^2) = \sum\limits_{k=1}^n(\sum\limits_{j=1}^n m_{kj}^2 m_{jk}^1) = tr (M_2 M_1).$ The lemma is proved.

{\bf Lemma 2.1.} {\it Let $H_j, \ph j=1,2$ be Hermitian matrices such that $H_j\ge 0, \ph j=1,2.$ Then
$$
tr (H_1 H_2) \ge 0.
$$
}

Proof. Let $U$ be an unitary transformation such that $\widetilde{H}_1 \equiv U H_1 U^* = diag \{h_1,\dots,h_n\}$.
Since any unitary transformation preserves the nonnegative definiteness of any Hermitian  matrix    we have
$$
h_j \ge 0, \phh j=\overline{1,n}. \eqno (2.1)
$$
Let $\widetilde{H}_2 \equiv U h_2 U^* = (h{ij})_{i,j=1}^n.$ As for as $\widetilde{H}_2$ is Hermitian it follows that (see [4], p. 300, Theorem 20)  $h_{jj} \ge 0, \ph j=\overline{1,n}$. This together with (2.1) implies
$$
tr (H_1 H_2) = tr([U H_1 U^*] [U H_2 U^*]) = tr (\widetilde{H}_1 \widetilde{h}_2) = \sum\limits_{j=1}^n h_j h_{jj} \ge 0.
$$
The Lemma is proved.

{\bf Lemma 2.3}. {\it Let $H$ be a Hermitian matrix of dimension $n \times n$ and let $K$ be a skew symmetric matrix ($K^* = - K$) of the same dimension. Then
$$
tr (H K) = 0. \eqno (2.3)
$$
}

Proof. By Lemma 2.1 we have
$$
tr(HK) = tr (K H) = tr ((K H)^*) = tr (H^* K^*) = - tr (H K).
$$
From here it follows (2.3). The lemma is proved.

{\bf Lemma 2.4.} {\it Let $H \ge 0$ be a Hermitian matrix of dimension $n \times n$ and let $V$ be any matrix of the same dimension. Then
$$
V H V^* \ge 0. \eqno (2.4)
$$
}

Proof. For any vectors $x$ and $y$ of dimension $n$ denote by $<x,y>$ their scalar product. Then $<V H V^* x, x> = <H(V^* x), (V^* x)> \ge 0$ (since $H \ge 0$). Hence (2.4) is valid. The lemma is proved.

Along with Eq. (1.1) consider the linear matrix system
$$
\sist{\Phi' = R(t) \Phi + P(t) \Psi,}{\Psi' = S(t)\Phi - Q(t) \Psi,} \ph t \ge t_0. \eqno (2.5)
$$
By a solution of this system we man a pair $(\Phi(t),\Psi(t))$ of absolutely continuous matrix functions of dimension $n\times n$ on $[t_0,+\infty)$, satisfying (2.5) almost everywhere on $[t_0,+\infty)$. It is not difficult to verify that all solutions of Eq. (1.1), existing on $[t_1,t_2) \subset [t_0,+\infty)$ are connected with solutions $(\Phi(t),\Psi(t))$ of the system (2.5) by relations
$$
\Phi'(t) = [R(t) + P(t)Y(t)]\Phi(t), \phh \Psi(t) = Y(t) \Phi(t), \phh t \in [t_1,t_2).  \eqno (2.6)
$$
By the Liouville's formula from here it follows
$$
\det \Phi(t) = \det \Phi(t_0)\exp\biggl\{\il{t_1}{t} tr (R(\tau) + P(\tau) Y(\tau)) d \tau\biggr\}, \phh t\in [t_1,t_2), \eqno (2.7)
$$
$$
\overline{\det \Phi(t)} = \overline{\det \Phi(t_0)}\exp\biggl\{\il{t_1}{t} tr (R^*(\tau) + Y^*(\tau)P(\tau) d \tau\biggr\}, \phh t\in [t_1,t_2). \eqno (2.8)
$$

{\bf 3. Global solvability criteria}.

\vskip 5pt
{\bf Definition 3.1.} {\it An interval $[t_1,t_2) \subset [t_0,+\infty)$ is called the maximum existence interval for a solution $Y(t)$ of Eq. (1.1), if $Y(t)$ exists on $[t_1,t_2)$ and cannot be continued to the right from $t_2$ as a solution of Eq. (1.1).
}

For any absolutely continuous matrix function $\Lambda(t)$ of dimension $n\times n$ on $[t_0,\infty)$  we set

$$
S_\Lambda(t) \equiv S(t) - \Lambda'(t) - \Lambda(t) P(t) \Lambda(t) - Q(t) \Lambda(t) - \Lambda(t) R(t), \phh t \ge t_0.
$$

{\bf Theorem 3.1.} {\it Let the following  conditions be satisfied.

\noindent
I) $P(t) \ge 0, \ph t \ge t_0.$

\noindent
For some absolutely continuous matrix function $\Lambda(t)$ of dimension $n\times n$ on $[t_0,\infty)$

\noindent
II)   $R(t) = Q^*(t) + P(t)(\Lambda^*(t) - \Lambda(t)) +  \mu(t) I, \ph t \ge t_0$, where $\mu(t)$ is a complex-valued locally integrable function on $[t_0,+\infty)$.

\noindent
III) $S_\Lambda (t) + S_\Lambda ^*(t) \ge 0, \ph t\ge t_0.$

\noindent
Then every solution $Y(t)$ of Eq. (1.1) with $Y(t_0) + Y^*(t_0) \ge \Lambda(t_0) + \Lambda^*(t_0)$ exists on $[t_0,+\infty)$ and
$$
Y(t) + Y^*(t) \ge \Lambda(t) + \Lambda^*(t), \phh t \ge t_0. \eqno (3.1)
$$
}

Proof.  In Eq. (1.1) substitute
$$
Y = Z + \Lambda(t), \phh t \ge t_0. \eqno (3.2)
$$
We obtain
$$
Z' + Z P(t) Z + Q_\Lambda(t) Z + Z R_\Lambda(t) - S_\Lambda(t) = 0, \phh t \ge t_0, \eqno (3.3)
$$
where $Q_\Lambda(t)\equiv Q(t) + \Lambda(t) P(t), \phh R_\Lambda(t) \equiv R(t) + P(t) \Lambda(t), \phh t \ge t_0.$
Let $Y(t)$ be  a solution of Eq. (1.1) with $Y(t_0) + Y^*(t_0) > \Lambda(t_0) + \Lambda^*(t_0)$ and let $[t_0,t_1)$ be its maximum existence interval.   Then by (3.2) $Z(t)\equiv Y(t) + \Lambda(t)$ is a solution of Eq. (3.3) on  $[t_0,t_1)$ and  $[t_0,t_1)$  is its maximum existence interval. Show that
$$
Z(t) + Z^*(t) > 0, \phh t \in[t_0,t_1). \eqno (3.4)
$$
Suppose this is not so. Then there exists $t_2 \in (t_0,t_1)$ such that
$$
Z(t) + Z^*(t) > 0, \ph t \in [t_0,t_2), \eqno (3.5)
$$
and
$$
\det [Z(t_2) + Z^*(t_2)] = 0. \eqno (3.6)
$$
Since $Z(t)$ is a solution of Eq. (3.3) on $[t_0,t_1)$ we have
$$
Z'(t) + Z(t) P(t) Z(t) + Q_\lambda(t) Z(t) + Z(t) R_\Lambda(t) - S_\Lambda(t) = 0,
$$
$$
[Z^*(t)]' + Z^*(t) P(t) Z^*(t) + R_\Lambda^*(t) Z^*(t) + Z^*(t) Q_\Lambda^*(t) - S_\Lambda^*(t) = 0,
$$
almost everywhere on $[t_1,t_2)$. Summing up these equalities and making some arithmetic transformations we obtain
$$
[Z(t) + Z^*(t)]' + [Z(t) + Z^*(t)] P(t) [Z(t) + Z^*(t)] + \frac{[Q_\Lambda(t) + R_\Lambda^*(t)]}{2}[Z(t) + Z^*(t)] +
$$
$$
+[Z(t) + Z^*(t)] \frac{[Q_\Lambda^*(t) + R_\lambda(t)]}{2} +\frac{[Q_\Lambda(t) - R_\Lambda^*(t)]}{2}[Z(t) - Z^*(t)] +
$$
$$
+[Z(t) - Z^*(t)] \frac{[R_\Lambda(t) - Q_\Lambda^*(t)]}{2} -[S_\Lambda(t) + S_\Lambda^*(t)] - Z(t) P(t) Z^*(t) - Z^*(t) P(t) Z(t) = 0 \eqno (3.7)
$$
almost everywhere on $[t_0,t_1)$. It follows from (3.5) that $(Z(t) + Z^*(t))^{-1}$ exists on $[t_0,t_2)$. Then (3.7) allows  to interpret $U(t) \equiv Z(t) + Z^*(t), \ph t\in [t_0,t_2)$ as a solution of the following linear matrix differential equation
$$
U' + \biggl\{[Z(t) + Z^*(t)] P(t) + \frac{Q_\Lambda(t) + R_\Lambda^*(t)}{2} + [Z(t) + Z^*(t)]\frac{Q_\Lambda^*(t) + R_\Lambda(t)}{2}[Z(t) + Z^*(t)]^{-1} +
$$
$$
+\Bigl(\frac{Q_\Lambda(t) - R_\Lambda^*(t)}{2}[Z(t) - Z^*(t)] + [Z(t) - Z^*(t)]\frac{R_\Lambda(t) - Q_\Lambda^*(t)}{2}\Bigr) [Z(t) + Z^*(t)]^{-1}-
$$
$$
-\Bigl(S_\Lambda(t) + S_\Lambda^*(t) + Z(t) P(t) Z^*(t) + Z^*(t) P(t) Z(t)\Bigr) [Z(t) + Z^*(t)]^{-1}\biggr\} U = 0, \ph t \in [t_0,t_2).
$$
Then by the Liouville formula we have
$$
\det [Z(t) + Z^*(t)] = \det [Z(t_0) + Z^*(t_0)] \exp\biggl\{ -\il{t_0}{t} tr  \biggl([Z(\tau) + Z^*(\tau)] P(\tau) +
$$
$$
+ \frac{Q_\Lambda(\tau) + R_\Lambda^*(\tau)}{2} + [Z(\tau) + Z^*(\tau)]\frac{Q_\Lambda^*(\tau) + R_\Lambda(\tau)}{2}[Z(t) + Z^*(\tau)]^{-1}+
$$
$$
\Bigl(\frac{Q_\Lambda(\tau) - R_\Lambda^*(\tau)}{2}[Z(\tau) - Z^*(\tau)] + [Z(\tau) - Z^*(\tau)]\frac{R_\Lambda(\tau) - Q_\Lambda^*(\tau)}{2}\Bigr) [Z(\tau) + Z^*(\tau)]^{-1}-
$$
$$
-\Bigl(S_\Lambda(\tau) + S_\Lambda^*(\tau) + Z(\tau) P(\tau) Z^*(\tau) + Z^*(\tau) P(\tau) Z(\tau)\Bigr) [Z(\tau) + Z^*(\tau)]^{-1}\biggr)d\tau\biggr\} = 0, \eqno (3.8)
$$
$t\in [t_0,t_2)$. By the condition II) we have
$$
\frac{Q_\Lambda(t) - R_\lambda^*(t)}{2}(Z(t) - Z^*(t))   + (Z(t) - Z^*(t))\frac{R_\Lambda(t) - Q_\Lambda^*(t)}{2} = \frac{\mu(t) - \overline{\mu(t)}}{2}(Z(t) - Z^*(t)),
$$
$t \in [t_0,t_2)$.
Since $Z(t) - Z^*(t), \ph t \in[t_0,t_2)$ is skew symmetric and $(Z(t) + Z^*(t))^{-1}, \linebreak t\in [t_0,t_2)$ is Hermitian By Lemma 2.3 from the last equality it follows that
$$
tr \biggl\{\frac{Q_\Lambda(t) - R_\Lambda^*(t)}{2}(Z(t) - Z^*(t)) + (Z(t) - Z^*(t))\frac{R_\Lambda(t) - Q_\lambda^*(t)}{2}\biggr\}(Z(t) + Z^*(t))^{-1} = 0,  \eqno (3.9)
$$
$ t\in[t_0,t_2).$ By Lemma 2.4 from I) it follows that
$$
Z(t) P(t) Z^*(t) + Z^*(t) P(t) Z(t) \ge 0, \ph t \in [t_0,t_2].
$$
Then since $(Z(t) + Z^*(t))^{-1} > 0, \ph t \in [t_0,t_2)$ by Lemma 2.2 we have
$$
tr [(Z(t)P(t)Z^*(t) + Z^*(t)P(t)Z(t))(Z(t) + Z^*(t))^{-1}] \ge 0, \ph t \in [t_0,t_2). \eqno (3.10)
$$
By the same reason from the condition III) of the theorem we obtain
$$
tr[(S_\Lambda(t) + S_\lambda^*(t))(Z(t) + Z^*(t))^{-1}]\ge 0, \phh t \in [t_0,t_2). \eqno (3.11)
$$
Since $Z(t) + Z^*(t), \ph t \in [t_0,t_1)$ is absolutely continuous on $[t_0,t_2]$ we have
$$
\il{t_0}{t_2}tr [(Z(t) + Z^*(t))P(t)]d t \le c =const. \eqno (3.12)
$$
Obviously
$$
tr(Z(t) + Z^*(t))\frac{Q_\Lambda^*(t) + R_\Lambda(t)}{2}(Z(t) + Z^*(t))^{-1} = tr \frac{Q_\Lambda^*(t) + R_\Lambda(t)}{2}, \ph t\in [t_0,t_2).
$$
This together with (3.8)-(3.12) implies that
$$
\det (Z(t_2) + Z^*(t_2)) \ge \det (Z(t_0) + Z^*(t_0)) \exp\bigl\{c\} \ne 0,
$$
which contradicts (3.6). The obtained contradiction proves (3.4).
It follows from (3.4) that $Y(t)$ exists on $[t_0,t_1)$ and
$$
Y(t) + Y^*(t) > \Lambda(t) + \Lambda^*(t), \phh t \in [t_0,t_1). \eqno (3.13)
$$

Let $(\Phi(t),\Psi(t))$ be a solution of the system(2.5) with $\Phi(t_0) = I, \ph \Psi(t_0) = Y(t_0)$. Then by (2.7) and (2.8) we have
$$
|\det\Phi(t)|^2 = |\det\Phi(t)|^2 \exp\biggl\{\il{t_0}{t} tr \biggl[R(\tau) + R^*(\tau) + P(\tau)(Y(\tau) + Y^*(\tau))\biggr]d\tau\biggr\}, \phh t \in [t_0,t_1).
$$
By Lemma 2.2 from here and from the condition I) of the theorem it follows that
$$
|\det\Phi(t_1)|^2 = |\det\Phi(t_0)|^2 \exp\biggl\{\min\limits_{t\in[t_0,t_1]}\il{t_0}{t} tr \biggl[R(\tau) + R^*(\tau) \biggr]d\tau\biggr\}> 0.
$$
Hence, $\det \Phi(t) \ne 0, \ph t\in [t_0, t_1 +\varepsilon)$ for some $\varepsilon > 0$. By (2.6) from here it follows that $Y_1(t)\equiv \Psi(t)\Phi^{-1}(t), \ph t \in[t_0,t_1+\varepsilon)$ is a solution of Eq. (1.1) on $[t_0,t_1 +\varepsilon)$ and coincides with $Y(t)$ on $[t_0,t_1)$. By the uniqueness theorem it follows from here that $[t_0,t_1)$ is not the maximum existence interval for $Y(t)$. Therefore $t_1 = +\infty$. This together with (3.13) implies
$$
Y(t) + Y^*(t) > \Lambda(t) + \Lambda^*(t),  \phh t \ge t_0.
$$
Thus the theorem is proved in the particular case when $Y(t_0) + Y^*(t_0) > \Lambda(t_0) + \Lambda^*(t_0)$. Let us prove the theorem in the general case, when $Y(t_0) + Y^*(t_0) \ge  \Lambda(t_0) + \Lambda^*(t_0)$. For any $\varepsilon > 0$ denote by $Y_\varepsilon(t)$ the solution of Eq. (1.1), satisfying the initial condition
$$
Y\varepsilon(t_0) = Y(t_0) +\varepsilon I,
$$
where $Y(t)$ is a solution of Eq. (1.1), with $Y(t_0) + Y^*(t_0) \ge \Lambda(t_0) + \Lambda^*(t_0), \ph I$ is the $n\times n$ dimensional identity matrix. Obviously
$$
Y_\varepsilon(t_0) + Y_\varepsilon^*(t_0) \ge 2\varepsilon I + \Lambda(t_0) + \Lambda^*(t_0) >  \Lambda(t_0) + \Lambda^*(t_0).
$$
Then by already proven $Y_\varepsilon(t)$ exists on $[t_0,+\infty)$ and
$$
Y_\varepsilon(t) + Y_\varepsilon^*(t)   >  \Lambda(t) + \Lambda^*(t), \phh t \ge t_0. \eqno (3.14)
$$
Let $[t_0,T)$ be the maximum existence interval for $Y(t)$. Show that
$$
T = +\infty. \eqno (3.15)
$$
Suppose $T< +\infty$. Let $\lambda(t)$ be the least eigenvalue of $Y(t) + Y^*(t) - (\Lambda(t) + \Lambda^*(t))$ and $\lambda_\varepsilon(t)$ be the least eigenvalue of $Y_\varepsilon(t) + Y_\varepsilon^*(t)  - (\Lambda(t) + \Lambda^*(t)), \ph t \in [t_0,T)$. It follows from (3.14) that
$$
\lambda_\varepsilon(t) > 0, \phh t \ge t_0. \eqno (3.16)
$$
Since the solutions of Eq. (1.1) are continuously dependent on their initial values we have $\lambda_\varepsilon(t) \to \lambda (t)$ as $\varepsilon \to 0$. This together with (3.16) implies that $\lambda(t) \ge 0,\ph t \in [t_0,T).$ Therefore
$$
Y(t) + Y^*(t) \ge  \Lambda(t) + \Lambda^*(t), \phh t\in [t_0,T). \eqno (3.17)
$$
Let $(\Phi_0(t),\Psi_0(t))$ be a solution of the system(2.5) with $\Phi_0(t_0) = I, \ph \Psi_0(t_0) = Y(t_0)$. Then by (2.7) and (2.8) we have
$$
|\det\Phi_0(t)|^2 = |\det\Phi_0(t)|^2 \exp\biggl\{\il{t_0}{t} tr \biggl[R(\tau) + R^*(\tau) + P(\tau)(Y(\tau) + Y^*(\tau))\biggr]d\tau\biggr\}, \ph t \in [t_0,T).
$$
By Lemma 2.2 it follows from here and from  (3.17)
that
$$
|\det \Phi_0(T)|^2 \ge |\det \Phi_0(t_0)|^2\exp\biggl\{\min\limits_{t \in [t_0,T]}\il{t_0}{t} tr (R(\tau) + R^*(\tau))d\tau\biggr\} > 0.
$$
Hence, $\det \phi_0(t) \ne 0, \ph t \in [t_0,T+\delta)$ for some $\delta > 0$. Then by (2.6)  $Y_0(t) \equiv \Psi_0(t)\Phi_0^{-1}(t), \ph t \in [t_0,T+\delta)$ is a solution of Eq. (1.1) on $[t_0,T+\delta)$, which coincides with $Y(t)$ on $[t_0,T)$. By the uniqueness theorem it follows from here that $[t_0,T)$ is not the maximum existence interval for $Y(t)$, which contradicts our supposition. The obtained contradiction proves (3.15). From (3.15) and (3.17) it follows (3.1). The theorem is proved.

{\bf Remark 3.1.} {\it Theorem 3.1 is an extension of Theorem 1.1 in two directions:

\noindent
I)  the set of coefficients of Eq. (1.1) is extended;

\noindent
II) the set of initial values, for which the solutions of Eq. (1.1) are continuable to $+\infty$, is extended.
}

 If $P(t) > 0, \ph t \ge t_0$ and $P^{-1}(t)[Q^*(t) - R(t) + \mu(t) I]$ is an absolutely continuous skew symmetric matrix function on $[t_0,\infty)$, then one  can verify that the condition II) of Theorem 3.1 is fulfilled if we take $\Lambda(t) \equiv P^{-1}(t)\frac{Q^*(t) - R(t) + \mu(t) I}{2}, \ph t \ge t_0.$ Then from Theorem 3.1 we obtain

{\bf Corollary 3.1.} {\it Let  $P(t) > 0, \ph t \ge t_0$  and let for some complex-valued locally integrable function  $\mu(t)$ the matrix function $2\Lambda_0(t) \equiv P^{-1}(t)[Q^*(t) - R(t) + \mu(t) I], t \ge t_0$  be skew symmetric. If $S_{\Lambda_0}(t) + S_{\Lambda_0}^*(t) \ge 0, \ph t \ge t_0,$ then every solution $Y(t)$ of the system (1.1) with $Y(t_0) + Y^*(t_0) \ge 0$ exists on $[t_0,\infty)$ and
$$
Y(t) + Y^*(t) \ge 0, \phh t \ge t_0.
$$
}

Proof. Note that $\Lambda_0(t) + \Lambda_0^*(t) \equiv 0, \ph t \ge t_0.$ The corollary is proved.

Let $\sqrt{P(t)}$ be absolutely continuous on $[t_0,+\infty)$. Multiply both sides of Eq. (1.1) at left and at right by $\sqrt{P(t)}$. Taking into account the equality
$$
\sqrt{P(t)}Y'\sqrt{P(t)} = (\sqrt{P(t)}Y-\sqrt{P(t)})' - \sqrt{P(t)}'Y\sqrt{P(t)} - \sqrt{P(t)}Y\sqrt{P(t)}'
$$
we obtain
$$
(\sqrt{P(t)}Y-\sqrt{P(t)})' + (\sqrt{P(t)}Y-\sqrt{P(t)})^2 + [\sqrt{P(t)} Q(t) - \sqrt{P(t)}'] Y \sqrt{P(t)} +
$$
$$
+ \sqrt{P(t)} Y [R(t) \sqrt{P(t)} - \sqrt{P(t)}'] - \sqrt{P(t)} S(t) \sqrt{P(t)} = 0, \phh t \ge t_0. \eqno (3.18)
$$
Consider the linear matrix equations
$$
\sqrt{P(t)} Q(t) - \sqrt{P(t)}' = X \sqrt{P(t)}, \phh t \ge t_0, \eqno (3.19)
$$
$$
R(t) \sqrt{P(t)} - \sqrt{P(t)}' = \sqrt{P(t)} X, \phh t \ge t_0. \eqno (3.20)
$$
Let $F(t)$ and $L(t)$ be absolutely continuous solutions of Eq. (3.19) and Eq. (3.20) respectively. Then from (3.18) we obtain
$$
Z' + Z^2 + F(t)Z + Z L(t) - \sqrt{P(t)}S(t)\sqrt{P(t)} = 0, \phh t \ge t_0, \eqno (3.21)
$$
where $Z\equiv \sqrt{P(t)} Y \sqrt{P(t)}$.

{\bf Corollary 3.2.} {\it Let  $P(t) > 0, \ph t \ge t_0, \ph  \sqrt{P(t)}, \ph t \ge t_0$ be absolutely continuous  and let for some complex-valued locally integrable on $[t_0,\infty)$ function $\nu(t)$ the matrix function
$$
T_\nu(t)\equiv \frac{\sqrt{P(t)}^{-1}[Q^*(t) - R(t)]\sqrt{P(t)} + \nu(t) I}{2}
$$
be absolutely continuous and skew symmetric. If
$$
\sqrt{P(t)}(S(t) + S^*(t))\sqrt{P(t)}  +  2 T_\nu^2(t) + (\overline{\nu}(t) - \nu(t))T_\nu(t)  \ge 0, \phh t\ge t_0. \eqno (3.22)
$$
then  every solution $Y(t)$ of Eq. (1.1) with
$$
\sqrt{P(t_0)}[Y(t_0) + Y^*(t_0)]\sqrt{P(t)} \ge 0
$$
exists on $[t_0,+\infty)$ and
$$
\sqrt{P(t)}[Y(t_0) + Y^*(t)]\sqrt{P(t)} \ge 0, \phh t \ge t_0.
$$
}

Proof.
If $P(t) > 0, \ph t \ge t_0$, then for the unique solutions $F(t)$ and $L(t)$ of the equations (3.19) and (3.20) respectively we have the equalities
$$
F(t) = [\sqrt{P(t)} Q(t) - \sqrt{P(t)}']\sqrt{P(t)}^{-1}, \phh t \ge t_0, \eqno (3.23)
$$
$$
L(t) = \sqrt{P(t)}^{-1}[R(t)\sqrt{P(t)} - \sqrt{P(t)}'], \phh t\ge t_0. \eqno (3.24)
$$
Set
$$
D_\nu(t)\equiv T_\nu'(t) + T_\nu^2(t) + [\sqrt{P(t)} Q(t) - \sqrt{P(t)}']\sqrt{P(t)}^{-1} T_\nu(t) +
$$
$$
+ T_\nu(t)\sqrt{P(t)}^{-1}[R(t)\sqrt{P(t)} - \sqrt{P(t)}'] - \sqrt{P(t)} S(t)\sqrt{P(t)}, \phh t \ge t_0.
$$
 In Eq. (3.21) substitute
$$
Z = U + T_\nu(t). \eqno (3.25)
$$
We obtain
$$
U' + U^2 + (F(t) + T_\nu(t)) U + U(L(t) + T_\nu(t)) + D_\nu(t) = 0, \phh t \ge t_0. \eqno (3.26)
$$
Since by the condition of the theorem  $T_\nu(t)$ is skew symmetric
$$
L(t) + T_\nu(t) = F^*(t) + T_\nu^*(t) + \nu(t) I, \phh t \ge t_0.
$$
Indeed, if  $T_\nu(t)$ is skew symmetric then the last equality is equivalent to
$$
2T_\nu(t) = F^*(t) - L(t) + \nu(t) I, \phh t \ge t_0,
$$
which is a straightforward consequence of (3.23) and (3.24). Hence, Eq. (2.6) takes the form
$$
U' + U^2 + (F(t) + T_\nu(t)) U + U(F^*(t) + T^*_\nu(t) + \nu(t) I) + D_\nu(t) = 0, \phh t \ge t_0. \eqno (3.27)
$$
By Theorem 3.1 every solution $U(t)$ of this equation with $U(t_0) + U^*(t_0) \ge 0$ exists on $[t_0,+\infty)$ and   $U(t) + U^*(t) \ge 0, \ph t \ge t_0$ provided
$$
D_\nu(t) + D_\nu^*(t) \le 0, \ph t \ge t_0. \eqno (3.28)
$$
Them in virtue of (3.25) every  solution $Z(t)$ of Eq. (3.21) with $Z(t_0) + Z^*(t_0) = U(t_0) + U^*(t_0) + T_\nu(t_0) + T_\nu^*(t_0) = U(t_0) + U^*(t_0) \ge 0$ (since $T_\nu(t_0) + T_\nu^*(t_0) = 0$) and
$$
Z(t) + Z^*(t) = U(t) + U^*(t) + T_\nu(t) + T_\nu^*(t) = U(t) + U^*(t) \ge 0, \phh t \ge t_0,
$$
provided (3.29) holds. Therefore due the connection $Z = \sqrt{P(t)} Y \sqrt{P(t)}$ between  (3.21) and (1.1) every solution $Y(t)$ of Eq. (1.1) with  $\sqrt{P(t_0)}[Y(t_0) + Y^*(t_0)]\sqrt{P(t)} \ge 0$ exists on $[t_0,+\infty)$ and  $\sqrt{P(t)}[Y(t_0) + Y^*(t)]\sqrt{P(t)} \ge 0, \phh t \ge t_0,$ provided (3.28) holds. Thus to complete the proof of the corollary  it remains to show that the conditions (3.28) and (3.22) are equivalent.  We have
$$
D_\nu(t) + D_\nu^*(t) = T_\nu'(t) + [T_\nu^{*}(t)]' + 2 T^2(t) + [\sqrt{P(t)} Q(t) - \sqrt{P(t)}']\sqrt{P(t)}^{-1} T_\nu(t) +
$$
$$
+ T_\nu(t)\sqrt{P(t)}^{-1}[R(t)\sqrt{P(t)} - \sqrt{P(t)}'] - \sqrt{P(t)} S(t) \sqrt{P(t)} -
$$
$$
- T_\nu(t)\sqrt{P(t)}^{-1}[Q^*(t)\sqrt{P(t)} - \sqrt{P(t)}']  - [\sqrt{P(t)} R^*(t) - \sqrt{P(t)}']\sqrt{P(t)}^{-1} T_\nu(t) -
$$
$$
- \sqrt{P(t)} S^*(t) \sqrt{P(t)} = 2 T_\nu^2(t) + [\sqrt{P(t)}(Q(t) - R^*(t))\sqrt{P(t)}^{-1} + \overline{\nu}(t) I] T_\nu(t) -
$$
$$
- \overline{\nu}T_\nu(t) + T_\nu(t)[\sqrt{P(t)}^{-1}[R(t) - Q^*(t)]\sqrt{P(t)} -\nu(t) I] + \nu(t) I -
$$
$$
- \sqrt{P(t)}(S(t) + S^*(t))\sqrt{P(t)} = - 2 T_\nu^2(t) + (\nu(t) - \overline{\nu}(t))T_\nu(t) - \sqrt{P(t)}(S(t) + S^*(t))\sqrt{P(t)},
$$
$t \ge t_0$. Therefore, the conditions (3.28) and (3.22) are equivalent. The corollary is proved.

\vskip 20 pt

\centerline{ \bf References}

\vskip 20pt

\noindent
1. G. J. Butler, Ch. R. Johnson, H. Wolkowicz. Nonnegative solutions of a quadratic matrix  \linebreak \phantom{aa}  equation, arising from comparison theorems in ordinary differential equations, SIAM \linebreak \phantom{aa}   J. Algebraic Discrete Methods 6 (1993) 9--20.

\noindent
2. G. Freiling, A survey of nonsymmetric Riccati equations. Linear Algebra Appl.\linebreak \phantom{aa} 351-352 (2002) 243--270.

\noindent
3. G. Freiling, G. Jank, A. Sarychev, Non-blow-up conditions for Riccati type matrix  \linebreak \phantom{aa} differential difference equations. results Math. 37 (1998) 84--103.

\noindent
4. F. G.  Gantmacher,  Theory of matrices. Moskow, Nauka, 1966, 576 pages.

\noindent
5. M. Gevert, Some properties of solutions of Riccati matrix differential equations. \linebreak \phantom{aa} Fast Math. 24 (1994) 42--53.

\noindent
6. G. A. Grigorian, Two comparison criteria for scalar Riccati equations and their \linebreak \phantom{aa} applications.
Izv. Vyssh. Uchebn. Zaved. Mat., 2012, Number 11, 20–35.

\noindent
7. G. A. Grigorian,    Criteria of global solvability for Riccati scalar equations

\noindent
8. G. A. Grigorian,    Properties of solutions of the scalar Riccati equation with complex \linebreak \phantom{aa} coefficients and some their applications. Differ. Equ. Appl. 10 (2018), no. 3, 277–298.

\noindent
9. G. AGrigorian, On the reducibility of systems of two linear first-order  ordinary\linebreak \phantom{aa}  differential equations. Monatsh. Math. 195 (2021), no. 1, 107–117.

\noindent
10 G. A. Grigorian,  New reducibility criteria for systems of two linear first-order ordinary \linebreak \phantom{aa} differential equations. Monatsh. Math. 198 (2022), no. 2, 311–322.

\noindent
11. J. Juang, Global existence for an abstract Riccati initial-value problem with possibly \linebreak \phantom{aa} unbounded coefficients. J. Math. Anal. Appl. 166(1) (1992) 103--111.

\noindent
12. J. Juang, M. T. Lee, Comparison theorems for the matrix Riccati equation, Linear \linebreak \phantom{aa} Algebra Appl. 196 (1994) 183--191.

\noindent
13.   Knobloch H. W., Pohl M. On Riccati  matrix differential equations. results Math. \linebreak \phantom{aa}31 (1997) 337--364.

\noindent
14.    W. T. Reid, Riccati Differential Equations, Academic Press, New York, 1972.

\end{document}